\documentclass[11pt]{article}
\usepackage[utf8]{inputenc}

\date{}

\usepackage{amsmath,amsfonts,amssymb,amsthm}
\usepackage{mathrsfs}
\usepackage{latexsym}
\usepackage[english]{babel}
\usepackage[usenames,dvipsnames]{xcolor}
\usepackage{comment}
\usepackage{appendix}

\usepackage{graphicx}

\usepackage[T1]{fontenc}
\usepackage[utf8]{inputenc}

\usepackage{comment}

\usepackage{indentfirst}
\usepackage[top=1.25in,left=1.25in,right=1.25in]{geometry}
\newlength{\numone}
\newlength{\widone}
\newlength{\numtwo}
\newlength{\widtwo}
\newcounter{countp}

\newtheorem{thm}{Theorem}%[section]

\theoremstyle{definition}

\newtheorem{que}[thm]{Question}

\numberwithin{equation}{section}

\author{\Large{Riccardo W. Maffucci and Niels Willems}}
\newcommand{\Addresses}{{
		\footnotesize
		
		R.W.~Maffucci, \textsc{EPFL MA SB,
			Lausanne, Switzerland 1015}\par\nopagebreak\vspace{-0.35cm}
		\textit{E-mail address}, R.W.~Maffucci: \texttt{riccardo.maffucci@epfl.ch}
		
		N.~Willems, \textsc{EPFL MA SB,
			Lausanne, Switzerland 1015}\par\nopagebreak\vspace{-0.35cm}
		\textit{E-mail address}, N.~Willems: \texttt{niels.willems@epfl.ch}
		}}

\title{\Large{\uppercase{\bf On smallest $3$-polytopes of given graph radius}}}
\begin{document}
%\titleformat{\section}
%  {\Large\scshape}{\thesection}{1em}{}
%\titleformat{\subsection}
%  {\large\scshape}{\thesubsection}{1em}{}
\maketitle

\begin{abstract}
The $3$-polytopes are planar, $3$-connected graphs. A classical question is, for $r\geq 3$, is the $2(r-1)$-gonal prism $K_2\times C_{2(r-1)}$ the unique $3$-polytope of graph radius $r$ and smallest size? Under some extra assumptions, we answer this question in the positive.
\end{abstract}
{\bf Keywords:} Graph radius, Extremal problem, Planar graph, $3$-polytope, Prism.
\\
{\bf MSC(2010):} 05C35, 05C07, 05C75, 52B05, 52B10, 05C10.

\section{Introduction}
The $1$-skeletons of polyhedra are the $3$-polytopes or polyhedral graphs. In 1934, Steinitz and Rademacher \cite{radste} showed that these are the $3$-connected, planar graphs. In 1961, Tutte \cite{tutte1961theory} proved that if $G$ is a $3$-polytope of size (i.e. number of edges) $q$ that is not a pyramid, then either $G$ or its dual can be derived by adding an edge to a $3$-polytopal graph of size $q-1$. While constructing the first few $3$-polytopes via this algorithm, one notes that the unique smallest one (fewest edges) of graph radius $3$ is the cube, of radius $4$ the hexagonal prism, of radius $5$ the octagonal prism. One naturally asks the following.

\begin{que}
\label{q:1}
For $r\geq 3$, is the $2(r-1)$-gonal prism $K_2\times C_{2(r-1)}$ the unique $3$-polytope of graph radius $r$ and smallest size?
\end{que}

This seems to be a hard problem, and actually a version of it has been addressed previously. In 1981, Harant-Walther \cite{harant1981radius} showed that every 3-connected graph of radius $r$ and order $p$ verifies $r<p/4+O(\log(p))$ as $p\to\infty$. Moreover, they give the example of the $2(r-1)$-gonal prism as a lower bound \cite[Figure 1]{harant1981radius}.

The Harant-Walther upper bound (i.e $|E|\leq 6(r-1)$) has since been improved \cite{harant1993upper,inoue1996radius}.
%Later Harant \cite{harant1993upper} lowered the upper bound to $r\leq\frac{p}{4}+8$, and Inoue \cite{inoue1996radius} improved it to $r\leq \frac{p}{4}+3$.
%Notice that the articles we cited never mention planarity, which could potentially be a condition to improve the bound in our original problem.
In other related literature, Klee \cite{klee1976classification} constructed all $3$-connected, $3$-regular graphs of given odd \textit{graph diameter} and smallest size. He described these minimal graphs as an efficient way to arrange the stations of a communicating network such that if any two stations are incapacitated, the network is still connected, and such that, in case of a breakdown, each station can rely on exactly three others.

Notice that none of these works assume planarity. In the present paper, we assume planarity together with some other conditions on the graph, and show that under these assumptions, $K_2\times C_{2(r-1)}$ is indeed the only solution to Question \ref{q:1}. Intuitively, one expects solutions to Question \ref{q:1} to be 3-regular graphs, however this is not so clear, as a solution may have fewer vertices than the $2(r-1)$-gonal prism together with vertices of degree $\geq 4$. Even under a $3$-regularity assumption, answering Question \ref{q:1} without further hypotheses is probably difficult.

\begin{thm}
\label{thm:1}
Let $G$ be a $3$-polytope of graph radius $r\geq 3$ and minimal size, and $u\in V(G)$ a vertex of eccentricity $r$. We define
\[V_i:=\{v\in V(G): d(u,v)=i\}, \qquad\qquad 0\leq i\leq r.\]
Suppose that
\begin{align}
&|V_r|=1;\label{h:1}\\ 
&|V_1|=|V_{r-1}|=3;\label{h:2}\\
 &|V_i|\geq 4&\text{for every } 2\leq i \leq r-2\label{h:3}.
\end{align}
Then $G\simeq K_2\times C_{2(r-1)}$.
\end{thm}

Section \ref{sec:pf} is dedicated to proving Theorem \ref{thm:1}. It can actually be unconditionally shown that, for each $2\leq i\leq r-2$, at least one of $|V_{i-1}|,|V_i|,|V_{i+1}|$ is $\geq 4$ \cite{harant1993upper,inoue1996radius}.

\paragraph{Notation and conventions.}
The term `graph' refers to an undirected, simple graph. If $G=(V,E)$ denotes a graph, the order of $G$ is $|V|$ and the size of $G$ is $|E|$. $K_n$ denotes the complete graph on $n$ vertices while $C_n$ refers to the cycle graph of order $n$. Let $G_1=(V_1,E_1)$ and $G_2=(V_2,E_2)$. %The join graph $G=G_1+G_2$ is defined as follows: $V_G=V_1\cup V_2$ and $E_G=E_1\cup E_2 \cup \{uv:u\in V_1, v\in V_2\}$.  The wheel graph on $n$ vertices, denoted $W_n$, is a graph isomorphic to $K_1+C_{n-1}$.
 %\begin{figure}[!h]
  %  \centering
    %\includegraphics[height=3cm]{wheel}
   % \caption{$W_7=K_1+C_6$, the wheel graph on 7 vertices.}
    %\label{fig:my_label}
%\end{figure}
 The Cartesian product $H=G_1\times G_2$ satisfies $V_H=V_1\times V_2$ and $E_H$ is defined as follows: if $(u_1,u_2),(v_1,v_2)\in V_H$, then $(u_1,u_2)(v_1,v_2) \in E_H$ if either $u_1=v_1$ and $u_2v_2\in E_2$, or $u_2=v_2$ and $u_1v_1\in E_1$.% \\
 %\begin{figure}[!h]
  %   \centering
     %\includegraphics[height=3cm]{cartesian}
   %  \caption{$K_2\times K_3$.}
    % \label{fig:my_label}
% \end{figure}
 %\\
\\A graph is said to be connected if for every pair of vertices $x\neq y$, we may find a path (i.e. a succession of incident edges) starting from $x$ and ending in $y$. 
 A graph $G$ of order at least $k+1$ is said to be $k$-connected if for every subset of vertices  $S$ of cardinality $k-1$ or less, $G-S$ is connected.\\
 %Intuitively, we say that a graph is planar if you can draw it on a sheet of paper without having 2 edges crossing. For example, $K_4$ is planar while $K_5$ is not (Figure 3).
%\begin{figure}[!h]
%    \centering
    %\includegraphics[height=5cm]{planar}
%    \caption{$K_4$ $(left)$ is planar while $K_5$ $(right)$ is not, as the edge $UV$ has to cross another edge}
 %   \label{fig:my_label}
%\end{figure}
 A graph is said to be $k$-regular if all vertices have degree $k$. % (that is, every vertex is adjacent to exactly $k$ vertices).
 A graph is said to be polyhedral if it is both planar and 3-connected. Notice that a 3-connected graph cannot have any vertex of degree $\leq 2$. %, as removing the neighbours of this vertex would make the obtained graph disconnected.
 The radius of a graph is the minimal eccentricity over all vertices of the graph, that is, if $G=(V,E)$, $\text{rad}(G)=\text{min}\{\text{ecc}_G(x) : x\in V\}$, with $\text{ecc}_G(x)=$max$\{d(x,y) : y\in V\}$.
The diameter of a graph is the maximal eccentricity or equivalently, the largest distance between two vertices of the graph.

\paragraph{Acknowledgements.}
R.M. was supported by Swiss National Science Foundation project 200021\_184927, held by Prof. Maryna Viazovska.
\\
N.W. worked on this project as partial fulfilment of his bachelor semester project at EPFL, spring 2022, under the supervision of R.M.

\section{Proof of Theorem \ref{thm:1}}
\label{sec:pf}
\subsection{Preliminaries}
Let $G=(V,E)$ satisfy the hypotheses of Theorem \ref{thm:1}. Our first remark is that \[|V|=\sum_{i=0}^r |V_i|=1+3+4(r-3)+3+1=4(r-1),\]
so that $G$ matches the order of the $3$-regular $K_2\times C_{2(r-1)}$. Thereby, as we want to minimise $|E(G)|$, we cannot have any vertex of degree $\geq 4$ in $G$, i.e. $G$ is 3-regular. Similarly, under our hypotheses, assumption \eqref{h:3} is equivalent to $|V_i|=4$ for every $2\leq i\leq r-2$. We also note that if $x\in V_j$, $y\in V_k$, and $xy\in E$, then $|k-j|\leq 1$.

%\subsection{}
Using \eqref{h:1}, call $w$ the only vertex at distance $r$ from $u$. By Menger's Theorem \cite{mccuaig1984simple}, a graph $G$ of order at least $k+1$ is $k$-connected if and only if between every pair of vertices there are $k$ pairwise internally disjoint paths. Since $G$ is $3$-connected, there are $3$ internally disjoint paths between $u$ and $w$. These $3$ paths are of length exactly $r$: they cannot be of length $<r$ as $d(u,w)=r$, and if one of them were of length $>r$, then there would be some other vertex on this path that is at distance $r$ from $u$, contradicting \eqref{h:1}. Call
\[x_1,x_2,...,x_{r-1},\]
\[y_1, y_2,... y_{r-1},\]
and
\[z_1,z_2,...,z_{r-1}\]
the vertices $\neq u,w$ along these three paths, where each $x_i,y_i,z_i\in V_i$.

By assumption, there exists $v\in V$ s.t. $d(x_1,v)=r$, and necessarily $v\in V_{r-1}$. By \eqref{h:2}, $v=y_{r-1}$ or $v=z_{r-1}$. Analogous statements may be made for $y_1,z_1$. This leaves two options up to relabelling: either
\[d(x_1,z_{r-1})=d(y_1,x_{r-1})=d(z_1,y_{r-1})=r,\]
or
\[d(x_1,z_{r-1})=d(y_1,x_{r-1})=d(z_1,x_{r-1})=r.\]

\subsection{First scenario}
\begin{figure}[h!]
    \centering
    \includegraphics[height=3.5cm]{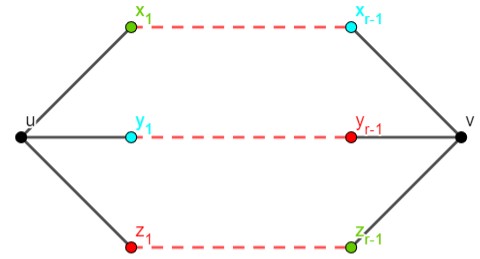}
    \caption{Vertices of the same colour are at maximal distance, $d(x_1,z_{r-1})=d(y_1,x_{r-1})=(z_1,y_{r-1})=r$. }
    \label{fig:3col}
\end{figure}

Here we have
\begin{equation}
\label{green}
d(x_1,z_{r-1})=r,
\end{equation}
\begin{equation}
\label{red}
d(y_1,x_{r-1})=r,
\end{equation}
and
\begin{equation}
\label{black}
d(z_1,y_{r-1})=r.
\end{equation}
This situation is depicted in Figure \ref{fig:3col}.

%Hence, 3 colours are necessary, but they are sufficient, and for the rest of the report, we will assume w.l.o.g. that the coloring is as in Figure 7:\\

Firstly, for all $1\leq i\leq r-1$, we have $x_iz_i\not\in E$ due to \eqref{green}, $y_ix_i\not\in E$ due to \eqref{red}, and $z_iy_i\not\in E$ due to \eqref{black}. Similarly, for all $1\leq i\leq r-2$, we have $x_iz_{i+1},y_ix_{i+1},z_iy_{i+1}\not\in E$. 

\begin{figure}[h!]
    \centering
    \includegraphics[height=3.5cm]{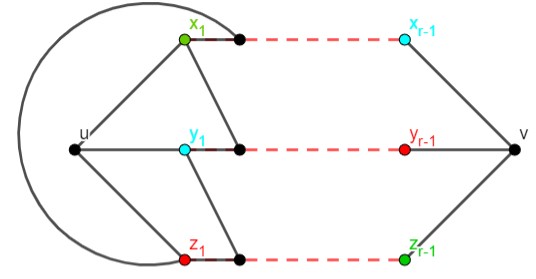}
    \caption{This configuration implies $|V_2|=3$.}
    \label{fig:biz}
\end{figure}

Next, we note that having all of $x_1y_2, y_1z_2, z_1x_2\in E$ at the same time is impossible (Figure \ref{fig:biz}): indeed, by definition each $v\in V_2$ is adjacent to at least one vertex in $V_1$, and since $G$ is $3$-regular, then $|V_2|=3$, contradicting \eqref{h:3}.
Moreover, by planarity one cannot have $v\in V_2$ such that $x_1v,y_1v,z_1v\in E$, and similarly one cannot have $v\in V_{r-2}$ such that $vx_{r-1},vy_{r-1},vz_{r-1}\in E$. Up to relabelling, this leaves two possibilities.

\begin{itemize}
\item
We have \[x_1y_2, y_1z_2, z_1t_2\in E\] for $t_2\in V_2$, $t_2\neq x_2,y_2,z_2$, and moreover \[y_{r-1}x_{r-2}, x_{r-1}t_{r-2}, z_{r-1}t_{r-2}\in E\] for $t_{r-2}\in V_{r-2}$, $t_{r-2}\neq x_{r-2},y_{r-2},z_{r-2}$. We will denote by $t_i$ the fourth vertex in $V_i$, $t_i\neq x_i,y_i,z_i$, $2\leq i\leq r-2$.

Now $t_{2}$ has two more neighbours, but $t_{2}x_2,t_{2}x_3,t_{2}y_3\not\in E$ by \eqref{black}. Then we have both of $t_{2}z_3,t_{2}t_3\in E$. Similarly, $t_{3}$ has two more neighbours, but $t_{3}x_3,t_{2}x_4,t_{2}y_3,t_{2}y_4\not\in E$ by \eqref{black}, and $t_{3}x_2\not\in E$ by \eqref{green}. Also, $t_{3}z_2,t_{3}z_3\not\in E$ as $z_2,z_3$ already have three neighbours. Then we have both of $t_{3}z_4,t_{3}t_4\in E$. We continue in this fashion, to see that
\[t_{i}z_{i+1},t_{i}t_{i+1}\in E, \qquad\qquad 2\leq i\leq r-3.\]
By $3$-regularity and the above considerations, we now deduce that
\[x_iy_{i+1}\in E,\qquad\qquad 2\leq i\leq r-3.\]
We have indeed constructed a $2(r-1)$-gonal prism in this case (Figure \ref{fig:my_label}).
\begin{figure}[h!]
    \centering
    \includegraphics[height=5.5cm]{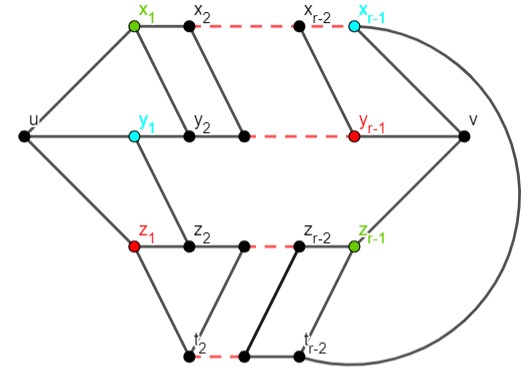}
    \caption{A $2(r-1)$-gonal prism.}
    \label{fig:my_label}
\end{figure}

\item
We have \[x_1y_2, y_1t_2, z_1t_2\in E.\]
Now $z_{r-1}x_{r-2},z_{r-1}y_{r-2}\not\in E$ by \eqref{green}, hence
$z_{r-1}t_{r-2}\in E$.

Further, $x_{r-1}y_{r-2}\not\in E$ by \eqref{red}, $y_{r-1}z_{r-2}\not\in E$ by \eqref{black}, and $x_{r-1}t_{r-2}$ and $y_{r-1}t_{r-2}$ cannot both be edges due to planarity. We are left with three cases to inspect in turn.
\begin{itemize}
\item
Say that
\[x_{r-1}z_{r-2},y_{r-1}t_{r-2}\in E.\]
Then $t_{r-2}x_{r-3},t_{r-2}y_{r-3}\not\in E$ by \eqref{green}, and $t_{r-2}z_{r-3}\not\in E$ by \eqref{black}. This leaves $t_{r-2}t_{r-3}\in E$. Proceeding in the same fashion, we get
\[t_{i}t_{i-1}\in E\qquad\qquad i=r-2,r-3,\dots,2.\]
We have constructed a $z_1y_{r-1}$-path of length $r-2$, contradicting \eqref{black} -- cf. Figure \ref{fig:new}.
\begin{figure}[h!]
    \centering
    \includegraphics[height=4cm]{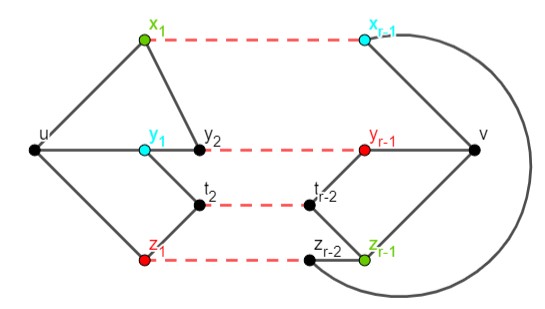}
    \caption{Here $d(z_1,y_{r-1})=r-2$.}
    \label{fig:new}
\end{figure}
\item
Say that
\[x_{r-1}t_{r-2},y_{r-1}x_{r-2}\in E.\]
This situation is also impossible: reasoning as above, one has to have a $t_{r-2}z_j$-path of length $r-2-j$ or $r-2-j+1$, and also a $t_2z_k$-path of length $k-2$, where $k\leq j-1$. This is because each $t_i$ is adjacent to a vertex in $V_{i-1}$ by definition. Altogether, there is a $y_1x_{r-1}$-path of length
\[\leq 1+(k-2)+(j-k)+(r-2-j+1)+1=r-1,\]
contradicting \eqref{red}.
\item
Say that
\[x_{r-1}z_{r-2},y_{r-1}x_{r-2}\in E.\]
By the above considerations, the third vertex adjacent to $z_i$ is $t_{i+1}$, to $t_i$ is $t_{i+1}$, and to $x_i$ is $y_{i+1}$, for each $i=2,3,\dots,r-3$. We have again obtained the $2(r-1)$-gonal prism.
\end{itemize}
\end{itemize}

\subsection{Second scenario}
Here we have
\begin{equation}
\label{colour}
d(x_1,z_{r-1})=d(y_1,x_{r-1})=d(z_1,x_{r-1})=r.
\end{equation}
Now $x_{r-1}$ is adjacent to an element of $V_{r-2}$ by definition. By \eqref{colour}, $x_{r-1}t_{r-2}\in E$. Next, either $t_{r-2}t_{r-3}\in E$ or $t_{r-1}x_{r-2}\in E$ (or both). In any case, we define
\begin{equation}
\label{min}
J:=\min\{j : x_jt_{j+1}\in E\}.
\end{equation}
Now $G$ is $3$-connected, hence $G-x_J-x_{r-1}$ is a connected graph. Therefore, there exists $k$ such that $t_kv\in E$, with $J+1\leq k\leq r-2$, and $v$ is one of \[y_{k-1},y_k,y_{k+1},z_{k-1},z_k,z_{k+1}.\]
By \eqref{colour}, $v=y_{k+1}$.

Still by $3$-connectivity, we must have either $x_ky_{k+1}\in E$, for some $J+1\leq k\leq r-2$, or $x_{J+1}t_{J}\in E$. But $x_ky_{k+1}\in E$ contradicts planarity, and $x_{J+1}t_{J}\in E$ contradicts the minimality of $J$. In any case we reach a contradiction, so that \eqref{colour} is impossible. The proof of Theorem \ref{thm:1} is complete.

\bibliographystyle{abbrv}
\bibliography{biblio}

\Addresses
\end{document}